\documentclass[a4paper,leqno]{amsart}
\usepackage{amsmath, amssymb, amsthm}

\usepackage[mathscr]{eucal}

\newtheorem{thm}{\textbf{Theorem}}[section]
\newtheorem{prop}[thm]{\textbf{Proposition}}
\newtheorem{lem}[thm]{\textbf{Lemma}}

\theoremstyle{definition}
\newtheorem{defn}[thm]{{\rm Definition}}
\newtheorem{conj}[thm]{{\rm Conjecture}}

\newcommand{\olapla}{\overline\bigtriangleup}
\newcommand{\onabla}{\overline\nabla}

\newcommand{\lapla}{\bigtriangleup}

\newcommand{\p}{\phi}

\title{$k$-harmonic maps into a Riemannian manifold with constant sectional curvature}

\author{Shun Maeta}

\curraddr{Nakakuki 3-10-9 Oyama-shi Tochigi
Japan}
\email{shun.maeta@gmail.com}

\subjclass[2000]{primary 58E20, secondary 53C43}

\begin{document}




\begin{abstract}
 J. Eells and L. Lemaire introduced k-harmonic maps,
 and Wang Shaobo showed the first variational formula.
 When, k=2, it is called biharmonic maps (2-harmonic maps). There have been
 extensive studies in the area.
 In this paper, we consider the relationship between 
 biharmonic maps and k-harmonic maps,
 and show non-existence theorem of 3-harmonic maps.
 We also give the definition of k-harmonic submanifolds of 
 Euclidean spaces, and study k-harmonic curve in Euclidean spaces. 
 Futhermore, we give a conjecture for k-harmonic submanifolds of Euclidean spaces.
\end{abstract}

\maketitle


\vspace{10pt}
\begin{flushleft}
{\large {\bf Introduction}}
\end{flushleft}
Theory of harmonic maps has been applied into various fields in differential geometry.
 The harmonic maps between two Riemannian manifolds are
 critical maps of the energy functional $E(\p)=\frac{1}{2}\int_M\|d\p\|^2v_g$, for smooth maps $\p:M\rightarrow N$.
 
On the other hand, in 1981, J. Eells and L. Lemaire \cite{jell1} proposed the problem to consider the {\em $k$-harmonic maps}:
 they are critical maps of the functional 
 \begin{align*}
 E_{k}(\p)=\int_Me_k(\p)v_g,\ \ (k=1,2,\dotsm),
 \end{align*}
 where $e_k(\p)=\frac{1}{2}\|(d+d^*)^k\p\|^2$ for smooth maps $\p:M\rightarrow N$.
G.Y. Jiang \cite{jg1} studied the first and second variational formulas of the bi-energy $E_2$, 
and critical maps of $E_2$ are called {\em biharmonic maps} ({\em 2-harmonic maps}). There have been extensive studies on biharmonic maps.
 
In 1989, Wang Shaobo \cite{ws1} studied the first variational formula of the 
$k$-energy $E_k$,
 whose critical maps are called $k$-harmonic maps.
 Harmonic maps are always $k$-harmonic maps by definition.
 But, the author \cite{sm1} showed biharmonic is not always $k$-harmonic $(k\geq 3)$.
 More generally, $s$-harmonic is not always $k$-harmonic $(s<k)$.
 Furthermore, the author \cite{sm1} showed the second variational formula of the $k$-energy.

In this paper, we study $k$-harmonic maps into a Riemannian manifold with constant sectional curvature $K$. 

In $\S \ref{preliminaries}$, we introduce notation and fundamental formulas of the tension field.

In $\S \ref{k-harmonic}$, we recall $k$-harmonic maps.

In $\S \ref{relationship}$, we give the relationship between biharmonic maps and $k$-harmonic maps.

In $\S \ref{3-harmonic}$, we study $3$-harmonic maps into a non positive sectional curvature and obtain non-existence theorem.

Finally, in $\S \ref{euclid}$, we define $k$-harmonic submanifolds of Euclidean spaces. And we show $k$-harmonic curve is a straight line.
 Furthermore, we give a conjecture for $k$-harmonic submanifolds in Euclidean spaces. 


\vspace{30pt}
\section{Preliminaries}\label{preliminaries}
Let $(M,g)$ be an $m$ dimensional Riemannian manifold,
 $(N,h)$ an $n$ dimensional one,
 and $\p:M\rightarrow N$, a smooth map.
 We use the following notation.
 The second fundamental form $B(\p)$
 of $\p$ is a covariant differentiation $\widetilde\nabla d\p$ of $1$-form $d\p$,
 which is a section of $\odot ^2T^*M\otimes \p^{-1}TN$.
For every $X,Y\in \Gamma (TM)$, let
 \begin{equation}
 \begin{split}
 B(X,Y)
&=(\widetilde\nabla d\p)(X,Y)=(\widetilde\nabla_X d\p)(Y)\\
&=\overline\nabla_Xd\p(Y)-d\p(\nabla_X Y)=\nabla^N_{d\p(X)}d\p(Y)-d\p(\nabla_XY). 
 \end{split}
 \end{equation}
 Here, $\nabla, \nabla^N, \overline \nabla, \widetilde \nabla$ are the induced connections on the bundles $TM$,
 $TN$, $\p^{-1}TN$ and $T^*M\otimes \p^{-1}TN$, respectively.
 
 If $M$ is compact,
 we consider critical maps of the energy functional
 \begin{align}
 E(\p)=\int_M e(\p) v_g,
 \end{align}
where $e(\p)=\frac{1}{2}\|d\p\|^2=\sum^m_{i=1}\frac{1}{2}\langle d\p(e_i),d\p(e_i)\rangle$
 which is called the {\em enegy density} of $\p$,
 the inner product 
 $\langle \cdot ,\cdot \rangle$ is a Riemannian metric $h$,
 and $\{e_i\}_{i=1}^m$ is a locally defined orthonormal frame field on $(M,g)$. 
 The {\em tension \ field} $\tau(\p)$ of $\p$ is defined by
 \begin{align}
 \tau(\p)=\sum^{m}_{i=1}(\widetilde \nabla d\p)(e_i,e_i)=\sum^m_{i=1}(\widetilde \nabla _{e_i}d\p)(e_i).
 \end{align}
 Then, $\p$ is a {\em harmonic map} if $\tau(\p)=0$.
 
 The curvature tensor field $ R^N(\cdot, \cdot)$ of the Riemannian metric on the bundle 
$TN$ is defined as follows :
 \begin{align}
R^N(X,Y)
=\nabla^N_X \nabla^N_Y - \nabla^N _Y \nabla^N_X-\nabla^N_{[X,Y]},
\ \ \ \ (X,Y\in \Gamma (TN)).
\end{align}
 

$\olapla
=\onabla^* \onabla
=-\sum^m_{k=1}(\onabla_{e_k}\onabla_{e_k}
-\onabla_{\nabla_{e_k}e_k}),$ 
is the {\em rough Laplacian}.



And G.Y.Jiang \cite{jg1} showed that $\phi:(M,g)\rightarrow (N,h)$ is a biharmonic (2-harmonic) if and only if 
$$\olapla \tau(\p) -R^N(\tau(\p),d\p(e_i))d\p(e_i)=0.$$



\vspace{30pt}

\section{$k$-harmonic maps}\label{k-harmonic}

J. Eells and L. Lemaire \cite{jell1} proposed the notation of $k$-harmonic maps. 
The Euler-Lagrange equation for the $k$-harmonic maps was shown by
Wang Shaobo \cite{ws1}.
In this section, we recall $k$-harmonic maps.

We consider a smooth variation $\{\p_{t}\}_{t\in I_{\epsilon}} (I_{\epsilon}=(-\epsilon, \epsilon))$ of $\p$ with parameter $t$,
 i.e., we consider the smooth map $F$  given by
$$F : I_{\epsilon}\times M\rightarrow N, F (t,p)=\p_{t}(p),$$
where $F (0,p)=\p_{0}(p)=\p(p), $ for all $p\in M$.

The corresponding variational vector field $V$ is given by
\begin{align*}
V(p)=\left.\frac{d}{dt}\right|_{t=0}\p_{t}(p)\in T_{\p(p)}N,
\end{align*}
$V$ are section of $\p^{-1}TN$, i.e., $V\in \Gamma (\p^{-1}TN)$.
 

\vspace{20pt}

\begin{defn}[\cite{jell1}]
For $k=1,2,\dotsm$ the {\em $k$-energy functional}
 is defined by 
\begin{align*}
E_k(\phi)=\frac{1}{2}\int_M\|(d+d^* )^k\phi\|^2v_g,\ \ \phi\in C^{\infty}(M,N).
\end{align*}
Then, $\phi$ is {\em $k$-harmonic} if it is a critical point of $E_k,$ i.e., for all smooth variation $\{\phi_t\}$ of $\phi$ with $\phi_0=\phi$,
\begin{align*}
\left.\frac{d}{dt}\right|_{t=0}E_k(\phi_t)=0.
\end{align*}
We say for a $k$-harmonic map to be {\em proper} if it is not harmonic. 
\end{defn}

\begin{thm}[\cite{ws1}]\label{2s-harmonic}
Let $k=2s\ (s=1,2,\cdots)$,
\begin{align*}
\left.\frac{d}{dt}\right|_{t=0}E_{2s}(\p_t)=-\int_M\langle \tau_{2s}(\p), V\rangle v_g,
\end{align*}
where, 
\begin{align*}
\tau_{2s}(\p)
=&\olapla^{2s-1}\tau(\p)
-R^N(\olapla^{2s-2}\tau(\p),d\p(e_j))e\p(e_j)\\
&-\sum^{s-1}_{l=1}
\{
R^N(\onabla_{e_j}\olapla^{s+l-2}\tau(\p),\olapla^{s-l-1}\tau(\p) )d\p(e_j) \\
&\hspace{32pt}-R^N(\olapla^{s+l-2}\tau(\p),\onabla_{e_j}\olapla^{s-l-1}\tau(\p) )d\p(e_j)
\},
\end{align*}
where, $\olapla^{-1}=0$, $\{e_i\}_{i=1}^m$ is a locally defined orthonormal frame field on $(M,g)$.
\end{thm}

\vspace{20pt}

\begin{thm}[\cite{ws1}]\label{2s+1-harmonic}
Let $k=2s+1\ \ \ (s=0,1,2,\cdots),$
\begin{align*}
\left.\frac{d}{dt}\right|_{t=0}E_{2s+1}(\p_t)=-\int_M\langle \tau_{2s+1}(\p), V\rangle v_g,
\end{align*}
where, 
\begin{align*}
\tau_{2s+1}(\p)
=&\olapla^{2s}\tau(\p)
-R^N(\olapla^{2s-1}\tau(\p),d\p(e_j))d\p(e_j)\\
&-\sum^{s-1}_{l=1}
\{
R^N(\onabla_{e_j}\olapla^{s+l-1}\tau(\p),\olapla^{s-l-1}\tau(\p) )d\p(e_j) \\
&\hspace{32pt}-R^N(\olapla^{s+l-1}\tau(\p),\onabla_{e_j}\olapla^{s-l-1}\tau(\p) )d\p(e_j)
\}\\
&-R^N(\onabla_{e_i}\olapla^{s-1}\tau(\p),\olapla^{s-1}\tau(\p))d\p(e_i),
\end{align*}
where, $\olapla^{-1}=0$, $\{e_i\}_{i=1}^m$ is a locally defined orthonormal frame field on $(M,g)$.
\end{thm}

\vspace{20pt}

\section{The relationship between biharmonic and $k$-harmonic}\label{relationship}

In \cite{sm1}, the auther showed $s$-harmonic is not always $k$-harmonic $(s<k)$.
Especially, biharmonic is not always $k$-harmonic $(k\geq 3)$.
So we study the relationship between biharmonic and $k$-harmonic $(2<k)$.
We obtain some results.

\begin{prop}\label{biharmonicii}
Let $\p: (M,g)\rightarrow (N,h)$ be an isometric immersion into a Riemannian manifold with constant sectional 
 curvature $K$. Then, $\p$ is biharmonic if and only if
$$\olapla\tau(\p)=Km\tau(\p).$$
\end{prop}

\vspace{10pt}

\begin{proof}
$\p$ is biharmonic if and only if 
\begin{align*}
0=&\olapla\tau(\p)-R^N(\tau(\p), d\p(e_i))d\p(e_i)\\
=&\olapla\tau(\p)-K\{\langle d\p(e_i),d\p(e_i) \rangle \tau(\p)-\langle d\p(e_i), \tau(\p)\rangle d\p(e_i)\}\\
=&\olapla\tau(\p)-Km\tau(\p).
\end{align*}
Thus, we have the proposition.
\end{proof}

\vspace{10pt}





\begin{lem}\label{innerl0}
Let $\p: (M,g)\rightarrow (N,h)$ be an isometric immersion into a Riemannian manifold with constant sectional 
 curvature $K$. If $\p$ is biharmonic, 
\begin{align}
\langle d\p(e_i), \olapla^l \tau(\p) \rangle =0.\ \ (l=0,1,\cdots )
\end{align}
\end{lem}

\vspace{10pt}

\begin{proof}
By using Proposition $\ref{biharmonicii}$, we have
\begin{align*}
\langle d\p(e_i), \olapla^l\tau(\p)\rangle 
=&mK\langle d\p(e_i),\olapla^{l-1}\tau(\p)\rangle \\
&\hspace{35pt}\cdots\\
=&(mK)^l\langle  d\p(e_i),\tau(\p)\rangle \\
=&0.
\end{align*}
\end{proof}

\vspace{10pt}

\begin{lem}\label{innernln}
Let $\p: (M,g)\rightarrow (N,h)$ be an isometric immersion into a Riemannian manifold with constant sectional 
 curvature $K$. If $\p$ is biharmonic, 
\begin{align}
\langle d\p(e_i),\onabla_{e_i}\olapla^l\tau(\p)\rangle =-(mK)^l||\tau(\p)||^2.
\end{align}
\end{lem}

\vspace{10pt}

\begin{proof}
By using Proposition $\ref{biharmonicii}$, we have
\begin{align*}
\langle d\p(e_i), \onabla_{e_i}\olapla^l \tau(\p)\rangle 
=&mK\langle d\p(e_i),\onabla_{e_i}\olapla^{l-1}\tau(\p)\rangle\\
&\hspace{50pt}\cdots\\
=&(mK)^l\langle d\p(e_i), \onabla_{e_i}\tau(\p)\rangle \\
=&-(mK)^l||\tau(\p)||^2,
\end{align*}

where, in the last equation, we only notice that
$$0=e_i\langle d\p(e_i),\tau(\p) \rangle =\langle \onabla_{e_i} d\p(e_i),\tau(\p)\rangle +\langle d\p(e_i),\onabla _{e_i}\tau(\p) \rangle .$$
\end{proof}

\vspace{10pt}

Using these lammas, we show the following two theorems.

\begin{thm}
Let $\p: (M,g)\rightarrow (N,h)$ be a biharmonic isometric immersion into a Riemannian manifold with constant sectional 
 curvature $K(\neq 0)$. If $\p$ is $2s$-harmonic $(s\geq 2)$, $\p$ is harmonic. 
\end{thm}

\vspace{10pt}

\begin{proof}
By Theorem \ref{2s-harmonic}, $\p$ is $2s$-harmonic if and only if 
\begin{align*}
&\olapla^{2s-1}\tau(\p)
-K\{m\olapla^{2s-2}\tau(\p)-\langle d\p(e_j),\olapla^{2s-2}\tau(\p) \rangle d\p(e_j)\}\\
&-\sum^{s-1}_{l=1}\{K(\langle \olapla^{s-l-1}\tau(\p),d\p(e_j)\rangle \onabla_{e_j}\olapla^{s+l-2}\tau(\p)\\
&\hspace{30pt}-\langle d\p(e_j),\onabla_{e_j}\olapla^{s+l-2}\tau(\p)\rangle \olapla^{s-l-1}\tau(\p)\\
&\hspace{30pt}-\langle \onabla_{e_j}\olapla^{s-l-1}\tau(\p),d\p(e_j)\rangle \olapla^{s+l-2}\tau(\p)\\
&\hspace{30pt}+\langle d\p(e_j),\olapla^{s+l-2}\tau(\p)\rangle \onabla_{e_j}\olapla^{s-l-1}\tau(\p))\}=0.
\end{align*}
By Proposition $\ref{biharmonicii}$, Lemma $\ref{innerl0}$ and $\ref{innernln}$, we have
\begin{align*}
0=&(mK)^{2s-1}\tau(\p)
-(mK)^{2s-1}\tau(\p)\\
&-\sum^{s-1}_{l=1}\{K((mK)^{2s-3}||\tau(\p)||^2\tau(\p)
+(mK)^{2s-3}||\tau(\p)||^2\tau(\p))\}\\
=&-2(s-1)K(mK)^{2s-3}||\tau(\p)||^2\tau(\p).
\end{align*}
Thus, we have the theorem.
\end{proof}

\vspace{10pt}

\begin{thm}
Let $\p: (M,g)\rightarrow (N,h)$ be a biharmonic isometric immersion into a Riemannian manifold with constant sectional 
 curvature $K(\neq 0)$. If $\p$ is $(2s+1)$-harmonic $(s\geq 1)$, $\p$ is harmonic. 
\end{thm}

\vspace{10pt}

\begin{proof}
By Theorem \ref{2s+1-harmonic}, $\p$ is $(2s+1)$-harmonic if and only if 
\begin{align*}
&\olapla^{2s}\tau(\p)
-K\{m\olapla^{2s-1}\tau(\p)-\langle d\p(e_j),\olapla^{2s-1}\tau(\p) \rangle d\p(e_j)\}\\
&-\sum^{s-1}_{l=1}\{K(\langle \olapla^{s-l-1}\tau(\p),d\p(e_j)\rangle \onabla_{e_j}\olapla^{s+l-1}\tau(\p)\\
&\hspace{30pt}-\langle d\p(e_j),\onabla_{e_j}\olapla^{s+l-1}\tau(\p)\rangle \olapla^{s-l-1}\tau(\p)\\
&\hspace{30pt}-\langle \onabla_{e_j}\olapla^{s-l-1}\tau(\p),d\p(e_j)\rangle \olapla^{s+l-1}\tau(\p)\\
&\hspace{30pt}+\langle d\p(e_j),\olapla^{s+l-1}\tau(\p)\rangle \onabla_{e_j}\olapla^{s-l-1}\tau(\p))\}\\
&-K\{\langle \olapla^{s-1}\tau(\p),d\p(e_i)\onabla_{e_i}\olapla^{s-1}\tau(\p)\\
&\hspace{30pt}-\langle d\p(e_i),\onabla_{e_i}\olapla^{s-1}\tau(\p)\rangle \olapla^{s-1}\tau(\p)\}=0.
\end{align*}
By Proposition $\ref{biharmonicii}$, Lemma $\ref{innerl0}$ and $\ref{innernln}$, we have
\begin{align*}
0=&(mK)^{2s}\tau(\p)
-(mK)^{2s}\tau(\p)\\
&-\sum^{s-1}_{l=1}\{K((mK)^{2s-2}||\tau(\p)||^2\tau(\p)+(mK)^{2s-2}||\tau(\p)||^2\tau(\p)\}\\
&-K(mK)^{2s-2}||\tau(\p)||^2\tau(\p)\\
=&-(2s-1)K(mK)^{2s-2}||\tau(\p)||^2\tau(\p).
\end{align*}
Thus, we have the theorem.
\end{proof}

\vspace{30pt}
\section{$3$-harmonic maps into non-positive curvature}\label{3-harmonic}
In this section we show non-existence theorem of $3$-harmonic maps.

G. Y. Jiang showed the follows.

\begin{thm}[\cite{jg1}]\label{nonposi biharmonic}
Assume that $M$ is compact and $N$ is non positive curvature,
 $i.e.,$ Riemannian curvature of $N$,  $K \leq 0.$
 Then, every biharmonic map $\p:M\rightarrow N$ is harmonic.
\end{thm}

We consider this theorem for $3$-harmonic maps.
First, we recall following theorem.

\begin{thm}[\cite{sm1}]\label{olapla harmonic}
Let $l=1,2,\dotsm $. If $\olapla^l\tau(\p)=0$ or $\onabla_{e_i}\olapla^{(l-1)}\tau(\p)=0,\ 
\newline(i=1,2,\dotsm,m)$, then
$\phi:M\rightarrow N$ from a compact Riemannian manifold into a Riemannian manifold is a harmonic map.
\end{thm}

\vspace{10pt}

Using this theorem, we obtain the next result.

\vspace{10pt}

\begin{prop}
Let $\p: (M,g)\rightarrow (N,h)$ be a isometric immersion from a compact Riemannian manifold into a Riemannian manifold with non positive constant sectional curvature $K\leq 0$. Then, $3$-harmonic is harmonic. 
\end{prop}

\vspace{10pt}

\begin{proof}
Indeed, by computing the Laplacian of the $4$-energy density $e_4(\p)$, we have 
\begin{equation}\label{le4}
\begin{split}
\lapla e_4(\p)
=&||\onabla_{e_i}\olapla\tau(\p)||^2-\langle \olapla^2\tau(\p), \olapla\tau(\p)\rangle \\
=&||\onabla_{e_i}\olapla\tau(\p)||^2 \\
&-\langle R^N(\olapla\tau(\p),d\p(e_i))d\p(e_i),\olapla\tau(\p)\rangle \\
&-\langle R^N(\onabla_{e_i},\tau(\p),\tau(\p))d\p(e_i),\olapla\tau(\p)\rangle, 
\end{split}
\end{equation}
due to $\p$ is $3$-harmonic.
Here, we consider the right hand side of 
$(\ref{le4})$, 
\begin{align*}
\langle R^N(\onabla_{e_i},\tau(\p),\tau(\p))d\p(e_i),\olapla\tau(\p)\rangle
=\langle &K\{\langle \tau(\p),d\p(e_i)\rangle \onabla_{e_i}\tau(\p)\\
&\hspace{5pt}-\langle d\p(e_i),\onabla_{e_i}\tau(\p)\rangle \tau(\p),\olapla\tau(\p)\}\rangle \\
&\hspace{-18pt}=K\{ ||\tau(\p)||^2\langle \tau(\p),\olapla\tau(\p)\rangle \}.
\end{align*}
Using Green's theorem, we have

\begin{equation}\label{nonnega}
\begin{split}
0=\int_{M}\lapla e_4(\p)
=&\int_M||\onabla_{e_i}\olapla\tau(\p)||^2\\
&-\langle R^N(\olapla \tau(\p),d\p(e_i))d\p(e_i),\olapla\tau(\p) \rangle\\
&-K||\tau(\p)||^2||\onabla_{e_j}\tau(\p)||^2v_g.
\end{split}
\end{equation}
Then, the both terms of $(\ref{nonnega})$ are non-negative, so we have

\begin{equation}
\begin{split}
0=\lapla e_4(\p)
=&||\onabla_{e_i}\olapla\tau(\p)||^2\\
&-\langle R^N(\olapla \tau(\p),d\p(e_i))d\p(e_i),\olapla\tau(\p) \rangle\\
&-K||\tau(\p)||^2||\onabla_{e_j}\tau(\p)||^2.
\end{split}
\end{equation}
Especially, we have
$$\onabla_{e_i}\olapla\tau(\p)=0.$$
Using Theorem $\ref{olapla harmonic}$, we obtain the proposition.
\end{proof}

\vspace{30pt}

\section{$k$-harmonic curves into Euclidean space}\label{euclid}
In this section, we consider $k$-harmonic curves into a Euclidean space $\mathbb{E}^n$ and we give a conjecture.
B. Y. Chen \cite{byc1} define biharmonic submanifolds of Euclidean spaces.

\begin{defn}[\cite{byc1}]
Let $x:M\rightarrow \mathbb{E}^n$ be an isometric immersion into a Euclidean space. 
$x:M\rightarrow \mathbb{E}^n$ is called biharmonic submanifold if 
\begin{equation}
\lapla^2 x=0,\ \text{that\ is},\ \lapla H=0,
\end{equation}
where, $H=-\frac{1}{m}\lapla x$ is the mean curvature vetor of the isometric immersion $x$ and 
 $\lapla$ the Laplacian of $M$. 
\end{defn}

B. Y. Chen and S. Ishikawa \cite{ci1} proved that any biharmonic surface in $\mathbb{E}^3$ is minimal.
And Chen \cite{byc1} gave a conjecture .

\begin{conj}[\cite{byc1}]
The only biharmonic submanifolds in Euclidean spaces are the minimal ones.
\end{conj}
There are several results for this conjecture (\cite{tijihu1}, \cite{id1} and \cite{ylo1} etc).
However, the conjecture is still open.
I. Dimitric \cite{id1} considered a cureve case $(n=1)$, and obtained following theorem.

\begin{thm}[\cite{id1}]
Let $x: C\rightarrow \mathbb{E}^n$ be a smooth curve parametrized by arc length, with the mean curvature vector $H$
 satisfying $\lapla  H=0$, then the curve is a straight line, i.e., totally geodesic in $\mathbb{E}^n$.
\end{thm}

\vspace{10pt}

We generalize this throrem.
First, we define $k$-harmonic submanifolds in Euclidean spaces.

\vspace{10pt}

\begin{defn}
Let $x:M\rightarrow \mathbb{E}^n$ be an isometric immersion into a Euclidean space.
$x:M\rightarrow \mathbb{E}^n$ is called $k$-harmonic submanifold if 
\begin{align}
\lapla ^k x=0,\ \text{that\ is},\ \lapla ^{k-1} H=0 \ (k=1,2,\cdots), 
\end{align}
where, $H=-\frac{1}{m}\lapla x$ is the mean curvature vetor of the isometric immersion $x$ and 
 $\lapla$ the Laplacian of $M$. 
\end{defn}

\vspace{10pt}

We also consider a curve case $(n=1)$, and obtain following theorem.

\vspace{10pt}

\begin{thm}
Let $x: C\rightarrow \mathbb{E}^n$ be a smooth curve parametrized by arc length, with the mean curvature vector $H$
 satisfying $\lapla ^{k-1} H=0, \ (k=1,2,\cdots )$ , then the curve is a straight line, i.e., totally geodesic in $\mathbb{E}^n$.
\end{thm}

\vspace{10pt}

\begin{proof}
We have 
$0=\lapla ^{k-1} H=-\lapla^k x=(-1)^{k+1}\frac{d^{2k}}{ds^{2k}}x, \ k=1,2,\cdots.$
Hence $x$ has to be a $(2k-1)$-th power polynomial in $s$, 
$$x=\frac{1}{2k-1}a_{2k-1}s^{2k-1}+\frac{1}{2k-2}a_{2k-2}s^{2k-2}+\cdots +a_1 s+a_0,$$
where $a_i\ (i=0,1,\cdots,2k-1)$ are constant vectors. Since $s$ is the natural parameter we have

\begin{align*}
1=&\langle \frac{dx}{ds}, \frac{dx}{ds} \rangle
=\langle \sum^{2k-1}_{i=1} a_{i}s^{i-1},\sum^{2k-1}_{i=1} a_{i}s^{i-1}\rangle\\
=&|a_{2k-1}|^2s^{4k-4}
+2\langle a_{2k-1}, a_{2k-2}\rangle s^{4k-5}
+\{2\langle a_{2k-1}, a_{2k-3} \rangle +|a_{2k-2}|^2\}s^{4k-6}\\
&\hspace{150pt}\cdots\\
&+\{2\langle a_1, a_3 \rangle +|a_2|^2\}s^2
+2\langle a_1, a_2\rangle s
+|a_1|^2
\end{align*}
On the right hand side we have a polynomial in $s$, so we must have
$$a_{2k-1}=a_{2k-2}=a_{2k-3}=a_{2k-4}=\cdots =a_{3}=a_2=0,\ |a_1|^2=1.$$
In other words, $x(s)=a_1s+a_0$ with $|a_1|^2=1$, and therefore the curve is a straight line.
\end{proof}

\vspace{10pt}

\begin{conj}
The only $k$-harmonic submanifolds in Euclidean spaces are the minimal ones.
\end{conj}

\vspace{10pt}

Especially, when $k=2$, it is B. Y. Chen conjecture.



\end{document}